
\magnification=1200

\catcode`\À=\active \defÀ{\`A}    \catcode`\à=\active \defà{\`a} 
\catcode`\Â=\active \defÂ{\^A}    \catcode`\â=\active \defâ{\^a} 
\catcode`\Æ=\active \defÆ{\AE}    \catcode`\æ=\active \defæ{\ae}
\catcode`\Ç=\active \defÇ{\c C}   \catcode`\ç=\active \defç{\c c}
\catcode`\È=\active \defÈ{\`E}    \catcode`\è=\active \defè{\`e} 
\catcode`\É=\active \defÉ{\'E}    \catcode`\é=\active \defé{\'e} 
\catcode`\Ê=\active \defÊ{\^E}    \catcode`\ê=\active \defê{\^e} 
\catcode`\Ë=\active \defË{\"E}    \catcode`\ë=\active \defë{\"e} 
\catcode`\Î=\active \defÎ{\^I}    \catcode`\î=\active \defî{\^\i}
\catcode`\Ï=\active \defÏ{\"I}    \catcode`\ï=\active \defï{\"\i}
\catcode`\Ô=\active \defÔ{\^O}    \catcode`\ô=\active \defô{\^o} 
\catcode`\Ù=\active \defÙ{\`U}    \catcode`\ù=\active \defù{\`u} 
\catcode`\Û=\active \defÛ{\^U}    \catcode`\û=\active \defû{\^u} 
\catcode`\Ü=\active \defÜ{\"U}    \catcode`\ü=\active \defü{\"u} 

\catcode`\ =\active \def { }

\hsize=11.25cm    
\vsize=18cm       
\parindent=12pt   \parskip=5pt     

\hoffset=.5cm   
\voffset=.8cm   

\pretolerance=500 \tolerance=1000  \brokenpenalty=5000

\catcode`\@=11

\font\eightrm=cmr8         \font\eighti=cmmi8
\font\eightsy=cmsy8        \font\eightbf=cmbx8
\font\eighttt=cmtt8        \font\eightit=cmti8
\font\eightsl=cmsl8        \font\sixrm=cmr6
\font\sixi=cmmi6           \font\sixsy=cmsy6
\font\sixbf=cmbx6

\font\tengoth=eufm10 
\font\eightgoth=eufm8  
\font\sevengoth=eufm7      
\font\sixgoth=eufm6        \font\fivegoth=eufm5

\skewchar\eighti='177 \skewchar\sixi='177
\skewchar\eightsy='60 \skewchar\sixsy='60

\newfam\gothfam           \newfam\bboardfam

\def\tenpoint{
  \textfont0=\tenrm \scriptfont0=\sevenrm \scriptscriptfont0=\fiverm
  \def\rm{\fam\z@\tenrm} \textfont1=\teni \scriptfont1=\seveni
  \scriptscriptfont1=\fivei
  \def\oldstyle{\fam\@ne\teni}\let\old=\oldstyle \textfont2=\tensy
  \scriptfont2=\sevensy \scriptscriptfont2=\fivesy
  \textfont\gothfam=\tengoth \scriptfont\gothfam=\sevengoth
  \scriptscriptfont\gothfam=\fivegoth \def\goth{\fam\gothfam\tengoth}
  
  \textfont\itfam=\tenit
  \def\it{\fam\itfam\tenit}
  \textfont\slfam=\tensl
  \def\sl{\fam\slfam\tensl}
  \textfont\bffam=\tenbf \scriptfont\bffam=\sevenbf
  \scriptscriptfont\bffam=\fivebf
  \def\bf{\fam\bffam\tenbf}
  \textfont\ttfam=\tentt
  \def\tt{\fam\ttfam\tentt}
  \abovedisplayskip=12pt plus 3pt minus 9pt
  \belowdisplayskip=\abovedisplayskip
  \abovedisplayshortskip=0pt plus 3pt
  \belowdisplayshortskip=4pt plus 3pt 
  \smallskipamount=3pt plus 1pt minus 1pt
  \medskipamount=6pt plus 2pt minus 2pt
  \bigskipamount=12pt plus 4pt minus 4pt
  \normalbaselineskip=12pt
  \setbox\strutbox=\hbox{\vrule height8.5pt depth3.5pt width0pt}
  \let\bigf@nt=\tenrm       \let\smallf@nt=\sevenrm
  \normalbaselines\rm}

\def\eightpoint{
  \textfont0=\eightrm \scriptfont0=\sixrm \scriptscriptfont0=\fiverm
  \def\rm{\fam\z@\eightrm}
  \textfont1=\eighti  \scriptfont1=\sixi  \scriptscriptfont1=\fivei
  \def\oldstyle{\fam\@ne\eighti}\let\old=\oldstyle
  \textfont2=\eightsy \scriptfont2=\sixsy \scriptscriptfont2=\fivesy
  \textfont\gothfam=\eightgoth \scriptfont\gothfam=\sixgoth
  \scriptscriptfont\gothfam=\fivegoth
  \def\goth{\fam\gothfam\eightgoth}
  
  \textfont\itfam=\eightit
  \def\it{\fam\itfam\eightit}
  \textfont\slfam=\eightsl
  \def\sl{\fam\slfam\eightsl}
  \textfont\bffam=\eightbf \scriptfont\bffam=\sixbf
  \scriptscriptfont\bffam=\fivebf
  \def\bf{\fam\bffam\eightbf}
  \textfont\ttfam=\eighttt
  \def\tt{\fam\ttfam\eighttt}
  \abovedisplayskip=9pt plus 3pt minus 9pt
  \belowdisplayskip=\abovedisplayskip
  \abovedisplayshortskip=0pt plus 3pt
  \belowdisplayshortskip=3pt plus 3pt 
  \smallskipamount=2pt plus 1pt minus 1pt
  \medskipamount=4pt plus 2pt minus 1pt
  \bigskipamount=9pt plus 3pt minus 3pt
  \normalbaselineskip=9pt
  \setbox\strutbox=\hbox{\vrule height7pt depth2pt width0pt}
  \let\bigf@nt=\eightrm     \let\smallf@nt=\sixrm
  \normalbaselines\rm}

\tenpoint

\def\pc#1{\bigf@nt#1\smallf@nt}         \def\pd#1 {{\pc#1} }

\catcode`\;=\active
\def;{\relax\ifhmode\ifdim\lastskip>\z@\unskip\fi
\kern\fontdimen2  -1.2 \fontdimen3 \font\fi\string;}

\catcode`\:=\active
\def:{\relax\ifhmode\ifdim\lastskip>\z@\unskip\fi\penalty\@M\ \fi\string:}

\catcode`\!=\active
\def!{\relax\ifhmode\ifdim\lastskip>\z@
\unskip\fi\kern\fontdimen2  -1.1 \fontdimen3 \font\fi\string!}

\catcode`\?=\active
\def?{\relax\ifhmode\ifdim\lastskip>\z@
\unskip\fi\kern\fontdimen2  -1.1 \fontdimen3 \font\fi\string?}

\catcode`\«=\active 
\def«{\raise.4ex\hbox{%
 $\scriptscriptstyle\langle\!\langle$}}

\catcode`\»=\active 
\def»{\raise.4ex\hbox{%
 $\scriptscriptstyle\rangle\!\rangle$}}

\frenchspacing

\def\raggedbottom{\topskip 10pt plus 36pt\r@ggedbottomtrue}

\def\pointir{\unskip . --- \ignorespaces}

\def\Medbreak{\vskip-\lastskip\medbreak}

\long\def\th#1 #2\enonce#3\endth{
   \Medbreak\noindent
   {\pc#1} {#2\unskip}\pointir{\it #3}\smallskip}

\def\decale#1{\smallbreak\hskip 28pt\llap{#1}\kern 5pt}
\def\decaledecale#1{\smallbreak\hskip 34pt\llap{#1}\kern 5pt}
\def\puce{\smallbreak\hskip 6pt{$\scriptstyle\bullet$}\kern 5pt}

\def\eqalign#1{\null\,\vcenter{\openup\jot\m@th\ialign{
\strut\hfil$\displaystyle{##}$&$\displaystyle{{}##}$\hfil
&&\quad\strut\hfil$\displaystyle{##}$&$\displaystyle{{}##}$\hfil
\crcr#1\crcr}}\,}

\catcode`\@=12

\showboxbreadth=-1  \showboxdepth=-1

\mathcode`A="7041 \mathcode`B="7042 \mathcode`C="7043 \mathcode`D="7044
\mathcode`E="7045 \mathcode`F="7046 \mathcode`G="7047 \mathcode`H="7048
\mathcode`I="7049 \mathcode`J="704A \mathcode`K="704B \mathcode`L="704C
\mathcode`M="704D \mathcode`N="704E \mathcode`O="704F \mathcode`P="7050
\mathcode`Q="7051 \mathcode`R="7052 \mathcode`S="7053 \mathcode`T="7054
\mathcode`U="7055 \mathcode`V="7056 \mathcode`W="7057 \mathcode`X="7058
\mathcode`Y="7059 \mathcode`Z="705A

\def\qp{{\bf Q}_p}

\def\azeroo#1{A_0(#1)_0}

\def\Kbar{{\overline K}}

\def\XKbar{{X{}_{\Kbar}}}

\def\XKtilde{X_\Ktilde}

\def\Ketoile{K^{\times}}
\def\Ktilde{{\widetilde K}}
\def\ktilde{{\tilde k}}

\def\kprim{{k'}}

\def\Gal{\mathop{\rm Gal}\nolimits}
\def\Pic{\mathop{\rm Pic}\nolimits}

\def\Hom{\mathop{\rm Hom}\nolimits}

\def\A{\mathord{\bf A}}

\def\P{\mathord{\bf P}}
\def\Z{\mathord{\bf Z}}

\def\ogoth{{\goth o}}
\def\ogothetoile{{{\goth o}^\times}}

\def\hfl#1#2#3{\smash{\mathop{\hbox to#3{\rightarrowfill}}\limits
^{\textstyle#1}_{\textstyle#2}}}
\def\gfl#1#2#3{\smash{\mathop{\hbox to#3{\leftarrowfill}}\limits
^{\textstyle#1}_{\textstyle#2}}}

\def\qed{\raise -2pt\hbox{\vrule\vbox to 10pt{\hrule width 4pt
                 \vfill\hrule}\vrule}}

\def\phi{\varphi}

\newcount\refno 

\long\def\ref#1:#2<#3>{                                        
\global\advance\refno by1\par\noindent                              
\llap{[{\bf\number\refno}]\ }{#1} \pointir{\it #2} #3\goodbreak }

\def\citer#1(#2){[{\bf\number#1}\if#2\empty\relax\else,\ #2\fi]}

\newcount\numerodesection
\def\section#1{\bigbreak
 {\bf\number\numerodesection.\ \ #1}\nobreak\medskip
 \advance\numerodesection by1}

\newcount\numeroderemarque
\def\remarque{\advance\numeroderemarque by1\smallbreak
{\it Remarque\/}\ \number\numeroderemarque~:}

\newcount\formuleno
\def\numeroter{\global\advance\formuleno by1
 \leqno{(\number\formuleno)}}
\def\formule(#1){{\rm (\number#1)}}

\def\rk{\mathop{\rm rk}\nolimits}

\def\somme{\mathop{\smash{\raise 2pt\hbox{$\sum$}}}\limits}

\def\long{\mathop{\rm long}\nolimits}

\newbox\bibbox
\setbox\bibbox\vbox{\bigbreak
\centerline{{\pc BIBLIOGRAPHICAL REFERENCES}}
 
\ref{\pc COLLIOT}-{\pc TH{\'E}L{\`E}NE} (J-L.) et {\pc CORAY} (D. F.):
L'{\'e}quivalence rationnelle sur les points ferm{\'e}s des surfaces
rationnelles fibr{\'e}es en coniques,
<Compositio Math., {\bf 39} (3), 1979, p.~301--332>
\newcount\ctcoray  \global\ctcoray=\refno

\ref{\pc COOMBES} (K. R.) and {\pc MUDER} (D. J.):
Zero-cycles on del Pezzo surfaces over local fields, 
<Journal of Algebra {\bf 97} (1985), 438--460.>
\newcount\coomuzero \global\coomuzero=\refno

\ref{\pc CORAY} (D.) and {\pc TSFASMAN} (M. A.):
Arithmetic on singular del Pezzo surfaces,
<Proceedings of the London Math.\ Soc.\ (3) {\bf 75} (1988), 25--87.>
\newcount\coraytsfas \global\coraytsfas=\refno

\ref{\pc CORTI} (A.):
{Del Pezzo surfaces over Dedekind schemes},
<Ann.\ of Math.\ (2) {\bf 144} (1996) 3, 641--683.>
\newcount\corti  \global\corti=\refno

\ref{\pc DALAWAT} (C. S.):
Le groupe de Chow d'une surface de Châtelet sur un corps local,
<Indag.\ mathem.\ N.S. {\bf 11} (2) (2000), 173--185, {\tt math.AG/0302156}.>
\newcount\chatelet  \global\chatelet=\refno

\ref{\pc DALAWAT} (C. S.):
Le groupe de Chow d'une surface rationnelle sur un corps local,
<{\tt math.AG/0302157}.>
\newcount\chow \global\chow=\refno

} 

\centerline{\bf The Chow group of a del Pezzo surface}
\centerline{\bf over a local field}

\vskip4mm

\centerline{Chandan Singh {\pc DALAWAT}}

\vskip2cm

The purpose of this Note is to give complete details of the example
sketched in \citer\chow() to illustrate the method for calculating the
Chow group of a rational surface over a local field split by an
unramified extension of the base.  The example concerned a del Pezzo
surface ; the construction of a regular integral model and the
determination of the specialisation map are carried out here.  As this
construction has been referred to in the literature \citer\corti(), it
might not be entirely out of place to make it more widely available.

So let $p\neq2$ be a prime number, $K$ a finite extension of $\qp$,
$\ogoth$ the ring of integers of $K$ and $v:\Ketoile\rightarrow\Z$ the
normalised valuation of~$K$.  Let $d\in\ogothetoile$ be a unit which
is not a square and let
$\alpha,\beta,\gamma,\delta\in\Ketoile$ be invertible elements of the
field $K$ whose valuations satisfy $v(\alpha)=v(\beta)+1$,
$v(\delta)=v(\gamma)+2$.  The condition on $d$ (that it be a unit) is
essential to the method ; the conditions on the valuations of
$\alpha,\beta,\gamma,\delta$ are not.  They could be related in other
ways but the case we have chosen suffices to illustrate the method.

In the projective space $\P_{4,K}$ with coordinates
$\def\sep{\,\colon}\,r\sep s\sep t\sep u\sep w$, consider the 
surface defined by :   
$$
\cases{\eqalign{
&\alpha(dr^2-s^2)=\beta(t-u)(t+w),\cr
&\gamma(dr^2-t^2)=\,\delta(s+u)(s+w);\cr}}
\numeroter\newcount\modelxok \global\modelxok=\formuleno 
$$
it is a del Pezzo surface of degree~4, birational over $K(\sqrt{d})$
to the projective plane.  The extension $K(\sqrt{d})$ is {\it
unramified\/} since $d\,$ is a unit.  So the natural map
$\Pic\XKtilde\rightarrow\Pic\XKbar\,$ is an isomorphism --- where
$\Kbar$ is an algebraic closure of $K$ and $\Ktilde$ is the maximal
unramified extension of $K$ in $\Kbar$ --- and the main technical
hypothesis of \citer\chow() is satisfied.

In order to apply the method of \citer\chow(th.~4), we first need to
construct a regular projective $\ogoth$-scheme $X$ whose generic fiber
$X_K$ is the surface \formule(\modelxok) and the irreducible
components of whose special fibre are --- with their reduced structure
--- smooth over $k$, the residue field of $K$.  We will then need to
calculate the image of the resulting specialisation map
\citer\chow({(30)}), which amounts to computing the intersection
numbers \citer\chow({(1)}) of the curves in the special fibre $X_k$ with
the various irreducible components of $X_{\ktilde}$, where $\ktilde$ is
the residue field of $\Ktilde$.
This is carried out in the following ten steps :
\medskip
\halign{\hskip2cm\hfil{\it #)}.\ &#\hfil\cr
i& The na{\"\i}ve model $X_0$\cr
ii& The first blow-up $X_1$\cr
iii& The first special fibre $X_{1,k}$\cr
iv& Singularities of the transform $X_1$\cr
v& The second blow-up $X$\cr
vi& The second special fibre $X_k$\cr
vii& Regularity of the transform $X$\cr
viii& The evaluation map\cr
ix& The Picard groups\cr
x& The specialisation map\cr
}
\smallbreak
{\it i).--- The na{\"\i}ve model $X_0$}.  Let $\pi$ be a uniformising
parameter of $K$.  Multiplying by powers of $\pi$ and by units, one
sees that the surface~\formule(\modelxok) is the generic fiber of the
projective $\ogoth$-surface $X_0$ defined in $E_0=\P_{4,\ogoth}$ by :
$$
\cases{\eqalign{
\pi(dr^2-s^2)-\beta(t-u)(t+w)=0,\cr
\gamma(dr^2-t^2)-\pi^2(s+u)(s+w)=0.\cr}}
\numeroter\newcount\modelxo \global\modelxo=\formuleno 
$$
Let us determine the singularities of the special fibre $X_{0,k}$ of
$X_0$, which is given in $\P_{4,k}$ by the system :
$$
(t-u)(t+w)=0,\quad {\overline d}r^2-t^2=0,
\numeroter\newcount\fibspexo \global\fibspexo=\formuleno 
$$
where $\overline x\in k$ denotes the reduction modulo $\pi\ogoth$ of
$x\in\ogoth$.  It has two irreducible components, $A_0$ and $B_0$,
defined in $X_{0,k}$ by :
$$
A_0:\quad t+w=0,\qquad B_0:\quad t-u=0.
$$
The union of the singular loci of $A_0$ and $B_0$ is defined in
$X_{0,k}$ by the system $r=0$, $t=0$.  Since we are interested in
constructing an $\ogoth$-surface the irreducible components of whose
special fibre (with their reduced structure) are smooth over $k$, let
us blow up the closed subscheme of $X_0$ defined by
$$
\pi=0,\quad r=0,\quad t=0.
\numeroter\newcount\singaobo \global\singaobo=\formuleno 
$$
But before doing that, let us see if there are points of $X_{0,k}$,
not on \formule(\singaobo), which are singular on $X_0$ : it comes
down to checking which points of the intersection $A_0\cap B_0$, apart
from the point $r=0$, $t=0$, $u=0$, $w=0$ which is already on
\formule(\singaobo), are in fact singular on $X_0$.  One verifies that
at all these points, with the possible exception of
$$
\eqalign{
R_0:\quad&\pi=0,\; s+t=0,\; t-u=0,\; t+w=0,\cr
S_0:\quad&\pi=0,\; s-t=0,\; t-u=0,\; t+w=0,\cr
}
$$
$A_0$ and $B_0$ are Cartier divisors on $X_0$ ; as they are regular at
these points, so is $X_0$.  On the other hand, a calculation shows
that the two points $R_0$ and $S_0$ are singular on $X_0$.

\smallbreak
{\it ii).--- The first blow-up $X_1$}.  Let us first construct the ambient
space of which $X_0$ blown up at \formule(\singaobo) is going to be a
closed subscheme.  Let $E_1=\P({\cal O}\oplus{\cal O}(1)\oplus{\cal
O}(1))$ be the projective bundle over $E_0=\P_{4,\ogoth}$ whose fibres
are projective planes $\P_{2,\ogoth}$ with coordinates
$\def\sep{\,\colon}\,P\sep R\sep T$.  The blow-up of $X_0$ along
\formule(\singaobo) is by definition the closed subscheme $X_1$
defined in $E_1$ by the system \formule(\modelxo) and the system
$$
\cases{\eqalign{
\pi R-rP=0,\ \pi T-tP=0,\  rT-tR=0,\cr
\gamma(dR^2-T^2)-P^2(s+u)(s+w)=0.\cr
}}
\numeroter\newcount\modelxi \global\modelxi=\formuleno 
$$

\smallbreak
{\it iii).--- The first special fibre $X_{1,k}$}.  The special fibre
$X_{1,k}$ of $X_1$ is defined in $E_{1,k}$ by the system
\formule(\fibspexo) and the reduction modulo $\pi$ of the system
\formule(\modelxi). The four irreducible components of $X_{1,k}$ are
defined in $X_1$ by the systems
$$
\eqalign{
A_1:\quad&P=0,\ t+w=0,\qquad& C_1:\quad&\pi=0,\ t=0,\ u=0,\cr
B_1:\quad&P=0,\ t-u=0,\qquad& D_1:\quad&\pi=0,\ t=0,\ w=0.\cr
}
$$
It is to be noted that, $\kprim$ being the quadratic extension of $k$
obtained by adjoining a square root of ${\overline d}$, there is
natural morphism $A_{1,\kprim}\rightarrow\P_{2,\kprim}$ (coordinates
$\def\sep{\,\colon}\,r\sep s\sep u$) whose fibre at a given closed
point $x$ is the pair of points $P=0$, ${\overline d}R^2-T^2=0$ of the
projective plane $\P_{2,\kprim}$ (coordinates
$\def\sep{\,\colon}\,P\sep R\sep T$) which is the fibre of
$E_1\rightarrow E_0$ at $x$.  The irreducible component $A_1$ thus
consists of two conjugate projective planes which do not meet ; in
particular, it is smooth over $k$.  The same holds for $B_1$.  

As for $C_1$, one sees that the fibre of the projection
$C_1\rightarrow\P_{1,k}$ (coordinates $\def\sep{\,\colon}\,s\sep w)$
at a given point $(s:w)$ is the conic
$$
\overline\gamma(\overline d R^2-T^2)-P^2s(s+w)=0
$$
in the plane $\P_{2,k}$ fibre of $E_1\rightarrow E_0$ at the said
point.  Thus, $C_1$ is a conic bundle over $\P_1$ with two degenerate
fibres, namely those above $(s:w)=(0:1)$ and $(1:-1)$ ; in particular,
$C_1$ is smooth and absolutely irreducible.  The same holds for $D_1$,
which is a conic bundle over $\P_{1,k}$ (coordinates
$\def\sep{\,\colon}\,s\sep u)$ whose fibre at any given point $(s:u)$
is 
$$
\overline\gamma(\overline d R^2-T^2)-P^2(s+u)s=0.
$$

It should also be noted that the intersection of $A_1$ and $B_1$ is a
pair of conjugate lines which do not meet and that $C_1$ and $D_1$
meet in the conic
$$
\overline\gamma(\overline d R^2-T^2)-P^2s^2=0
\numeroter\newcount\fibcom \global\fibcom=\formuleno
$$
in the fibral $\P_{2,k}$ at the point of intersection of the two lines
above which $C_1$ and $D_1$ are conic bundles.

\smallbreak
{\it iv).--- Singularities of the transform $X_1$}.  Let us determine the
singularities of the strict transform $X_1$ (given by the systems
\formule(\modelxo), \formule(\modelxi)).  As the projection
$X_1\rightarrow X_0$ is an isomorphism in the complement of $C_1\cup
D_1$ in $X_1$, the inverse images of $R_0$ and $S_0$, namely
$$
\eqalign{
R_1:\quad&P=0,\; s+t=0,\; t-u=0,\; t+w=0,\cr
S_1:\quad&P=0,\; s-t=0,\; t-u=0,\; t+w=0,\cr
}
$$
are singular on $X_1$.  All other singularities of $X_1$, if any, must
lie on $C_1\cup D_1$.  

Now, at every point of $C_1\cup D_1$ not in their intersection
$C_1\cap D_1$, the divisors $A_1$, $B_1$, $C_1$, $D_1$ are Cartier ;
since they are also regular, so is $X_1$ at those points.  The trace
of $C_1\cap D_1$ on the open subscheme $Y: P\neq0,\ s\neq0$ of $X_1$
is defined by the ideal ${\cal I}$ generated by $\pi$, $u$ and $w$.
The $({\cal O}_{Y}/{\cal I})$-module ${\cal I}/{\cal I}^2$ is free of
rank~2 : one has the relation
$$
\pi=-\beta(t-u)(t+w)+\pi dr^2
$$
which shows that $\pi\in{\cal I}^2$.  As the trace in
question (the conic \formule(\fibcom) deprived of its closed point
$P=0$) is regular of codimension~2, it follows that $X_1$ is regular
at these points.  

Finally, it turns out that the closed point defined in $C_1\cap D_1$
by $P=0$ and in $X_1$ by
$$
M_1:\quad P=0,\; t=0,\; u=0,\; w=0
$$
is indeed singular on $X_1$ ; it is also the point where the four
irreducible components of $X_{1,k}$ meet.  Thus the only singularities
of $X_1$ are the three points $R_1$, $S_1$, $M_1$ ; their union is
defined by
$$
P=0,\; t-u=0,\; t+w=0,\; t(dr^2-s^2)=0.
\numeroter\newcount\troissing \global\troissing=\formuleno
$$
\smallbreak
{\it v).--- The second blow-up $X$}.  We first construct an ambient space
$E_2$ from $E_1$ in the same manner as $E_1$ was constructed out of
$E_0$ in {\it ii)} above, but this time the fibres are $\P_{3,\ogoth}$
(coordinates $\def\sep{\,\colon}\,V\sep U\sep W\sep Z$) of
``\thinspace weights\thinspace'' $(0:1:1:3)$ (just as they were of
weights $(0:1:1)$ in {\it ii)} above).  In $E_2$, the blow-up $X$ of
$X_1$ along \formule(\troissing) is defined by the systems
\formule(\modelxo), \formule(\modelxi) and the system
$$
\cases{
\rk\pmatrix{V & U & W & Z \cr 
P & t-u & t+w & t(dr^2-s^2)\cr}<2,\cr
\hfil\cr
\quad\quad\qquad VZ-\beta T.UW=0.\cr
}
$$
\smallbreak
{\it vi).--- The second special fibre $X_k$}.  The seven irreducible
components of the special fibre $X_{k}$ are the following :
\medskip
\halign{\qquad$#\;$:&\quad#\hfil\cr
A& the blow-up of $A_1$ at the closed points $R_1$, $S_1$, $M_1$,\cr
B& the blow-up of $B_1$ at the closed points $R_1$, $S_1$, $M_1$,\cr
C& the blow-up of $C_1$ at the closed point $M_1$,\cr
D& the blow-up of $D_1$ at the closed point $M_1$,\cr
R& the fibre above the closed point $R_1$,\cr
S& the fibre above the closed point $S_1$,\cr
M& the fibre above the closed point $M_1$.\cr
}\noindent
where fibres of the projection $X\rightarrow X_1$ are meant.  All
seven irreducible components (with their reduced structure) are smooth
over $k$.  The first six are of multiplicity~1 ; let us show that $M$
is of multiplicity~2.  The functions $T$ and $dr^2-s^2$ are invertible
on $M$ and $t-u$ is a local equation for $M$ in the open subscheme of
$X$ where $T$, $dr^2-s^2$ and $U$ are invertible.  In this open
subscheme, we have
$\displaystyle \pi={tP\over T}$, 
$\displaystyle t={Z\over (dr^2-s^2)U}\cdot(t-u)$ and
$\displaystyle P={V\over U}\cdot(t-u)$,  
from which we get the relation
$$
\pi={VZ\over T(dr^2-s^2)U^2}\cdot(t-u)^2,
$$
showing that $M$ is of multiplicity~2 in $X_k$.  The divisor
corresponding to the special fibre thus consists of
$$\def\\#1{\quad\setbox0\hbox{M} \hbox to\wd0{\hfil #1\hfil}}
\eqalign{&\\A\\B\\C\\D\\R\\S\\2M\cr
&\\2\\2\\1\\1\\2\\2\phantom{2}\\2\cr}
\numeroter\newcount\fibspe \global\fibspe=\formuleno
$$
where the second row indicates that the only absolutely irreducible
components are $C$ and $D$, each of the five others splits as two
irreducible components over $\ktilde$ --- the residue field of
$\Ktilde$.  These are in fact defined over the quadratic extension of
$k$ in which $\overline d\in k$ is a square.

The the mutual intersections of these seven components are easy to
compute and will be of use later.

\smallbreak
{\it vii).--- Regularity of the transform $X$}.  Let $D$ be the union of
the three divisors $R$, $S$, $M$ of $X$ ; its complement is regular
since it is isomorphic to $X_1$ deprived of the closed subscheme
\formule(\troissing), which we know to be regular (cf.~{\it iv)\/}
above).  Now the Cartier divisor $D$ (defined by $P=0$ when $V\neq0$,
by $t-u$ when $U\neq0$, by $t-w$ when $W\neq0$ and by $t(dr^2-t^2)$
when $Z\neq0$) is regular and hence so is $X$ at each of its points.
This shows that the scheme $X$ is everywhere regular.

\smallbreak
{\it viii).--- The evaluation map} \citer\chow(p.~7).  Let $F$ be the
free $\Z$-module on the irreducible components of $X_\ktilde$ and let
$(f_Y)_Y$ be the canonical basis --- indexed by the irreducible
components of $X_k$ --- of $\Hom_k(F,\Z)$ (where $\Hom_k$ stands for
$\Gal(\ktilde|k)$-homomorphisms).  We have the evaluation map
$\xi:\Hom_k(F,\Z)\rightarrow\Z$ which to a
$\Gal(\ktilde|k)$-homomorphism $F\rightarrow\Z$ associates its value
at the special fibre $X_\ktilde$ considered as an element of $F$
(cf.~\formule(\fibspe)).  From what we have said about the
multilplicities of the components and about their absolute
irreducibilty, one has
$$
\def\\#1{\quad\setbox0\hbox{M} \hbox to\wd0{\hfil#1\hfil}}
\eqalign{
Y=&\\A\\B\\C\\D\\R\\S\\M\phantom{.}\cr
\xi(f_Y)=&\\2\\2\\1\\1\\2\\2\\{4.}\cr} 
\numeroter\newcount\xivalues \global\xivalues=\formuleno
$$
\smallbreak
{\it ix).--- The Picard groups}.  In order to complete the
calculation of the Chow group of the $K$-surface \formule(\modelxo)
--- the generic fibre of\/ $X$ --- by the method of \citer\chow(th.~4),
we wish to determine the Picard groups of the seven irreducible
components of $X_k$.

--- ${\Pic A}$.  Let ${\goth l}_A$ be the inverse image
under $A\rightarrow A_1$ of a pair of conjugate lines in $A_1$ not
passing through the points $R_1$, $S_1$, $M_1$.  Let ${\goth d}_{R,A}$
be the inverse image of a pair of conjugate lines in $A_1$ passing
through $R_1$ but distinct from the pair $A_1\cap B_1$.  Similarly,
define ${\goth d}_{S,A}$ (resp.~${\goth d}_{M,A}$) to be the inverse
image of a pair of conjugate lines in $A_1$ passing through $S$ but
distinct from $A_1\cap B_1$ (resp. passing through $M$ but distinct
from $A_1\cap D_1$).  The group $\Pic A$ is generated by the classes
of the four divisors ${\goth l}_A$, ${\goth d}_{R,A}$, ${\goth
d}_{S,A}$, ${\goth d}_{M,A}$.

--- ${\Pic B}$.  In the preceding discussion, if we replace
$A$ by $B$, $B$ by $A$,\penalty-1000 $C$ by $D$ and $D$ by $C$, we get
four divisors ${\goth l}_B$, ${\goth d}_{R,B}$, ${\goth d}_{S,B}$,
${\goth d}_{M,B}$ whose classes generate the group $\Pic B$.

--- ${\Pic C}$.  Let ${\goth f}_C$ be the inverse image
under $C\rightarrow C_1$ of a fiber of $C_1\rightarrow\P_1$ which does
not pass through the point $M_1$.  The group $\Pic C$ is generated by
the classes of the three divisors ${\goth f}_C$, $C\cap B$, $C\cap M$.

--- ${\Pic D}$.  In an analogous fashion, the group $\Pic D$
is generated by the classes of the divisors ${\goth f}_D$, $D\cap A$,
$D\cap M$.

--- ${\Pic R}$, ${\Pic S}$, ${\Pic M}$.
The group $\Pic R$ is genrated by the divisors $R\cap A$ and $R\cap
B$.  The same is true for $\Pic S$ and for $\Pic M$.

\smallbreak
{\it x).--- The specialisation map} \citer\chow(p.~7).  We need to
determine the images of the specialisation maps $\Pic Y \rightarrow
\Hom_k(F,\Z)$ \citer\chow({(30)}), where $Y$ runs through the
irreducible components of $X_k$.  In other words, we have to determine
the subgroup of $\Hom_k(F,\Z)$ generated by the images of the various
divisors on $X_k$ which we have listed in {\it ix)\/} above.  

But consider for example the {\it divisor\/ $C\cap B$ on\/} $C$.  Its
{\it image in\/} $\Pic B$ is a linear combination of the classes of
${\goth l}_B$, ${\goth d}_{R,B}$, ${\goth d}_{S,B}$, ${\goth d}_{M,B}$
and hence its image in $\Hom_k(F,\Z)$ is contained in the image of
$\Pic B$.

Ruling out such cases, we are reduced to computing the images of the
ten divisors which have had the honour of being designated by a gothic
letter.  These images are expressed in terms of the canonical basis
$(f_Y)_Y$ of $\Hom_k(F,\Z)$, where $Y$ runs through the irreducible
components \formule(\fibspe) of the special fibre $X_k$ : 
{\catcode`\+=\active
\def+{\hbox{\phantom{$-$}}}
$$\def\\#1{{\goth #1}}
\bordermatrix{&\\l_A&\\d_{R,A}&\\d_{S,A}&\\d_{M,A}&\\f_D&
\\l_B&\\d_{R,B}&\\d_{S,B}&\\d_{M,B}&\\f_C\cr
f_A&-2&-1&-1&-2&+1&+1\cr
f_B&+1&  &  &  &  &-2&-1&-1&-2&+1\cr
f_C&  &  &  &  &  &+2&  &  &  &-2\cr
f_D&+2&  &  &  &-2\cr
f_R&  &+1&  &  &  &  &+1\cr
f_S&  &  &+1&  &  &  &  &+1\cr
f_M&  &  &  &+1&  &  &  &  &+1\cr
} 
$$}

Notice that the ``\thinspace degree\thinspace'' of each column
is~$0$ ({\it i.e.}\ $\xi$ vanishes at that element of $\Hom_k(F,\Z)$),
in accordance with \citer\chow(lemma 10) --- for example the degree of
the first column is
$$
\xi(-2f_A+f_B+2f_D)=(-2).2+1.2+2.1=0
$$ 
(cf.~\formule(\xivalues)).  Now the map $\gamma_0$ \citer\chow({(16)}) is
an isomorphism of the Chow group $\azeroo{X_K}$ of $0$-cycles of
degree~0 on $X_K$ with the torsion subgroup of the quotient of
$\Hom_k(F,\Z)$ by the subgroup generated by the columns of the
displayed matrix \citer\chow().  Hence :

\th PROPOSITION 
\enonce 
The Chow group $\azeroo{X_K}$ of\/ $0$-cycles of dgree\/ $0$ on the
surface \formule(\modelxo) is isomorphic to\/ $\Z/2\Z$.
\endth

The  result is in conformity with Coombes and Muder
\citer\coomuzero(th.~4.5), who used methods specific to del Pezzo
surfaces (split by an unramified extension of the base).

If we use Coray and Tsfasman \citer\coraytsfas(prop.~4.6) to write
\formule(\modelxo) as the Châtelet surface
$$
y^2-dz^2=x(x-e_1)(x-e_2),\quad 
e_1=\gamma-\pi^2+{\pi^3\over \beta},\ 
e_2={\gamma\pi\over\beta},
\numeroter\newcount\chatmodel \global\chatmodel=\formuleno
$$
the result agrees with a calculation of Colliot-Thélène
\citer\coraytsfas(prop.~4.7) (see also \citer\chatelet(prop.~1)) as
the Chow group is a birational invariant \citer\ctcoray(prop.~6.3).

\medbreak

Let us terminate with a word about the case when $d\in\Ketoile$ is not
a unit in \formule(\modelxok), \formule(\chatmodel) ; the extension
$K(\sqrt{d})$ of $K$ is then {\it ramified}.  No general method seems
to be available for calculating the Chow group of a rational surface
over a local field which is not split by an unramified extension of
the base.  But in this particular case, the Chow group has been
calculated in \citer\chatelet(prop.~2) as a function of the type of
bad reduction (cf.~\citer\chatelet(prop.~5)).

\unvbox\bibbox
\bigskip\bigskip
\leftline{Chandan Singh Dalawat}
\leftline{Harish-Chandra Research Institute}
\leftline{Chhatnag Road, Jhunsi}
\leftline{{\pc ALLAHABAD} 211\thinspace019 India}

\bye